\documentclass[11pt]{article}
\usepackage{amsfonts}
\usepackage{bbm}
\usepackage{mathrsfs}
\usepackage{amssymb}
\usepackage{url}

\title{Structural Properties of The Quantized Matrix Algebra $D_q(n)$\\ Established by Means of Gr\"obner-Shirshov Basis theory\thanks{Project supported by the National Natural Science Foundation of China (1186106).}\\
{\normalsize --- dedicated to Professor Andr\'e Leroy on the occasion of his retirement}}

\author{Lina Niu and Rabigul Tuniyaz\thanks{E-mail: rabigul802@hotmail.com}\\
{\small Department of Mathematics, School of  Science}\\
{\small Xinjiang Institute of Science and Technology}\\
{\small  Akesu, 843100, Xinjiang, China}\\}

\date{}

\begin{document}
\maketitle
\begin{center}
\begin{minipage}{120mm}
{\small {\bf Abstract.} Let $D_q(n)$ be the quantized matrix algebra introduced by Dipper and Donkin.  It is shown that some structural properties of $D_q(n)$ and their modules may be established and realized by means of Gr\"obner-Shirshov basis theory.}
\end{minipage}\end{center} {\parindent=0pt\par

}\vskip -.5truecm

\renewcommand{\thefootnote}{\fnsymbol{footnote}}
\let\footnote\relax

\def\NZ{\mathbb{N}}
\def\QED{\hfill{$\Box$}}
\def \r{\rightarrow}
\def\mapright#1#2{\smash{\mathop{\longrightarrow}\limits^{#1}_{#2}}}

\def\v5{\vskip .5truecm}
\def\OV#1{\overline {#1}}
\def\hang{\hangindent\parindent}
\def\textindent#1{\indent\llap{#1\enspace}\ignorespaces}
\def\item{\par\hang\textindent}

\def\LH{{\bf LH}}\def\LM{{\bf LM}}\def\LT{{\bf
LT}}\def\KX{K\langle X\rangle} \def\KS{K\langle X\rangle}
\def\B{{\cal B}} \def\LC{{\bf LC}} \def\G{{\cal G}} \def\FRAC#1#2{\displaystyle{\frac{#1}{#2}}}
\def\SUM^#1_#2{\displaystyle{\sum^{#1}_{#2}}}
\def\T{\widetilde}\def\KD{K\langle D\rangle}
\def\PRC{\prec_{d\textrm{\tiny -}lex}}

\section*{1. Introduction}

Let $K$ be a field of characteristic 0. The quantized matrix algebra $D_q(n)$,   introduced in [2], has been widely studied and generalized in different contexts, for instance, [3], [4], [5], [6], and [7]. In this paper, we show that some structural properties of $D_q(n)$ may be established and realized in a constructive-computational way. More precisely, after this introduction section, in Section 2 we show explicitly that the defining relations of $D_q(n)$ form a Gr\"obner-Shirshov basis, and that this Gr\"obner-Shirshov basis determines a PBW $K$-basis $\B$ of $D_q(n)$ which is in the standard form. In Section 3, we show, by constructing an appropriate monomial ordering $\prec$ on its PBW $K$-basis ${\cal B}$ , that $D_q(n)$ is a solvable polynomial algebra in the sense of [9], thereby $D_q(n)$ has an algorithmic Gr\"obner basis theory for  (two-sided, one-sided) ideals and modules. To demonstrate how  the main results of Section 2 and Section 3  may bring some perspective of establishing and realizing structural properties of $D_q(n)$ and their modules in a constructive-computational way, in Section 4 we specify several applications of the main results of Section 2 and Section 3.\par

Throughout this note, $K$ denotes a field of characteristic 0, $K^*=K-\{0\}$, and all $K$-algebras considered are associative with multiplicative identity 1. If $S$ is a nonempty subset of an algebra $A$, then we write $\langle S\rangle$ for the two-sided ideal of $A$ generated by $S$.\v5

\section*{2. The defining relations of $D_q(n)$ form a Gr\"obner-Shirshov basis}
In this section, all terminologies concerning Gr\"obner-Shirshov bases, such as composition, ambiguity, and normal word, etc., are referred to [1].\par
Let $K$ be a field of characteristic 0, $I(n)=\{(i,j)~|~i,j=1,2,\cdots,n\}$ with $n\ge 2$, and let $D_q(n)$ be the quantized matrix algebra with the set of $n^2$ generators
$d=\{ d_{ij}~|~(i,j)\in I(n)\}$, in the sense of [2], namely, $D_q(n)$ is the associative $K$-algebra generated by the  given $n^2$ generators  subject to the relations:
$$\begin{array}{ll}
d_{ij}d_{st}=qd_{st}d_{ij},&\hbox{if}~i>s~\hbox{and}~j\le t,\\
d_{ij}d_{st}=d_{st}d_{ij}+(q-1)d_{sj}d_{it},&\hbox{if}~i>s~\hbox{and}~j>t,\\
d_{ij}d_{ik}=d_{ik}d_{ij},&\hbox{for~all}~i,j,k, \end{array}$$
where $i,j,k,s,t=1,2,...,n$ and $q\in K^*$ is the quantum parameter.\par

Now, let $D=\{D_{ij}~|~(i,j)\in I(n)\}$, $\KD$  the free associative $K$-algebra generated by $D$, and let $S$ denote the set of defining relations of $D_q(n)$ in  $\KD$, that is, $S$ consists of elements
$$\begin{array}{ll}
(a)~f_{ijst}=D_{ij}D_{st}-qD_{st}D_{ij},&\hbox{if}~i>s~\hbox{and}~j\le t,\\
(b)~g_{ijst}=D_{ij}D_{st}-D_{st}D_{ij}-(q-1)D_{sj}D_{it},&\hbox{if}~i>s~\hbox{and}~j>t,\\
(c)~h_{ijik}=D_{ij}D_{ik}-D_{ik}D_{ij},&\hbox{for~all}~i,j,k.\end{array}$$
Then, $D_q(n)\cong\KD /\langle S\rangle$ as $K$-algebras, where
$\langle S\rangle$ denotes the (two-sided) ideal of $\KD$ generated by $S$, i.e., $D_q(n)$ is presented as a quotient of $\KD$. Our aim below is to show that $S$ forms a Gr\"obner-Shirshov basis with respect to a certain monomial ordering on $\KD$. To this end,
let us take the deg-lex ordering $\PRC$ (i.e., the {\it degree-preserving  lexicographic ordering}) on $D^*$, where $D*$ is the set of all mono words in letters of $D$, i.e., all words of finite length like $u=D_{ij}D_{kl}\cdots D_{st}$. More precisely, we first take the  lexicographic ordering $<_{lex}$ on $D^*$ which is the natural extension of the ordering on the set $D$ of generators of $\KD$: for $D_{ij}$, $D_{kl}\in D$,
$$D_{kl}<D_{ij}\Leftrightarrow\left\{\begin{array}{l}k<i,\\ \hbox{or}~k=i~\hbox{and}~l<j,
\end{array}\right.$$
and for two words $u=D_{k_1l_1}D_{k_2l_2}\cdots D_{k_sl_s}$, $v=D_{i_1j_1}D_{i_2j_2}\cdots
D_{i_tj_t}\in D^*$,
$$\begin{array}{rcl} u\prec_{lex} v&\Leftrightarrow&\hbox{there exists an}~m~\hbox{such that}\\
&{~}&D_{k_1l_1}=D_{i_1j_1}, D_{k_2l_2}=D_{i_2j_2},\ldots , D_{k_{m-1}l_{m-1}}=D_{i_{m-1}j_{m-1}}\\
&{~}&\hbox{but}~D_{k_ml_m}<D_{i_mj_m}\end{array}$$
(note that conventionally the empty word $1<D_{ij}$ for all $D_{ij}\in D$). For instance
$$D_{41}D_{21}D_{31}\prec_{lex}D_{42}D_{13}D_{43}\prec_{lex}D_{42}D_{23}D_{34}D_{41}.$$
And then,  by assigning each $D_{ij}$ the degree 1, $1\le i,j\le n$, and writing $|u|$ for the degree of a word $u\in D^*$, we take the deg-lex ordering $\PRC$ on the set $D^*$: for $u,v\in D^*$,
$$u\PRC v\Leftrightarrow\left\{\begin{array}{l}|u|<|v|,\\ \hbox{or}~|u|=|v|~\hbox{and}~u\prec_{lex}v.
\end{array}\right.$$
For instance, $$D_{24}D_{11}\PRC D_{32}D_{24}\PRC D_{32}D_{31}\PRC D_{11}D_{12}D_{13}.$$
It is straightforward to check that $\PRC$ is a monomial ordering on $\KD$, namely, $\PRC$ is a well-ordering and
$$u\PRC v~\hbox{implies}~wur\PRC wvr~\hbox{for all}~u, v, w, r\in D^*.$$
With this monomial ordering $\PRC$ in hand, we are ready to prove the following result.\v5

{\bf Theorem 2.1} With notation as fixed above, let $J=\langle S\rangle$ be the ideal of $D_q(n)$ generated by $S$. Then, with respect to the monomial  ordering  $\PRC$ on $\KD$, the set $S$ is a Gr\"obner-Shirshov basis of the ideal $J$, i.e., the defining relations of $D_q(n)$ form a Gr\"obner-Shirshov basis.\vskip 6pt

{\bf Proof} By  [1], it is sufficient to check that  all  compositions determined by elements in $S$ are trivial modulo $S$. In doing so, let us first fix two more notations. For an element  $f\in\KD$, we write $\OV{f}$ for the leading mono word of $f$ with respect to $\PRC$, i.e., if $f=\sum_{i=1}^s\lambda_iu_i$ with $\lambda_i\in K$, $u_i\in D^*$, such that $u_1\PRC u_2\PRC\cdots\PRC u_s$, then $\OV{f}=u_s$. Thus, the set $S$ of defining relations of $D_q(n)$  has the set of leading mono words
$$\OV{S}=\left\{\begin{array}{l} \OV{f}_{ijst}=D_{ij}D_{st},~\hbox{if}~s<i,~j\le t,\\
\OV{g}_{ijst}=D_{ij}D_{st},~\hbox{if}~s<i,~t<j,\\
\OV{h}_{ijik}=D_{ij}D_{ik},~\hbox{if}~k<j.\end{array}\right\}$$
(note that if $j<k$, then since $h_{ijik}=D_{ij}D_{ik}-D_{ik}D_{ij}$, we have $\OV{h}_{ijik}=D_{ik}D_{ij}$). Also let us write $(a\wedge b)$ for the composition determined by defining relations (a) and (b) in $S$. Similar notations are made for compositions of other pairs of defining relations in $S$. \par

By means of $\OV{S}$ above, we start by listing all possible ambiguities $w$ of compositions of intersections determined by elements in $S$, as follows:
$$\begin{array}{lll} (a\wedge a)&w=D_{ij}D_{st}D_{kl}&\hbox{if}~i>s>k,~j\le t\le l,\\
(a\wedge b)&w_1=D_{ij}D_{st}D_{kl}&\hbox{if}~i>s>k,~j\le t,~t>l,\\
{~}&w_2=D_{ij}D_{st}D_{kl}&\hbox{if}~i>s>k,~j\le t,~t\le l,\\
(a\wedge c)&w_1=D_{ij}D_{st}D_{sk}&\hbox{if}~i>s,~j\le t,~t>k,\\
{~}&w_2=D_{ij}D_{ik}D_{st}&\hbox{if}~i>s,~j>k,~k\le t,\\
(b\wedge b)&w=D_{ij}D_{st}D_{kl}&\hbox{if}~i>s>k,~j>t>l,\\
(b\wedge c)&w_1=D_{ij}D_{st}D_{sk}&\hbox{if}~i>s,~j>t>k,\\
{~}&w_2=D_{ij}D_{ik}D_{st}&\hbox{if}~i>s,~j>k>t,\\
(c\wedge c)&w=D_{ij}D_{ik}D_{it}&\hbox{if}~j>k>t.\end{array}$$
Instead of writing down all tedious verification processes, below we shall record  only the verification processes of four typical cases:
$$\begin{array}{lll}
(a\wedge b)&w_1=D_{ij}D_{st}D_{kl},&\hbox{if}~i>s>k,~j\le t,~t>l,\\
(a\wedge c)&w_1=D_{ij}D_{st}D_{sk},&\hbox{if}~i>s,~j\le t,~t>k,\\
(b\wedge b)&w=D_{ij}D_{st}D_{kl},&\hbox{if}~i>s>k,~j>t>l,\\
(b\wedge c)&w_1=D_{ij}D_{st}D_{sk},&\hbox{if}~i>s,~j>t>k,\end{array}$$
because other cases can be checked in a similar way (the interested reader may contact the
author directly in order to see other verification processes).\par

$\bullet$ The case $(a\wedge b)$ with $w_1=D_{ij}D_{st}D_{kl}$, where $i>s>k,~j\le t,~t>l$.
\par
Since in this case $w_1=D_{ij}D_{st}D_{kl}=\bar{f}_{ijst}D_{kl}=D_{ij}\bar{g}_{stkl}$ with
$f_{ijst}=D_{ij}D_{st}-qD_{st}D_{ij}$, where $i>s$ and $j\le t$, and $g_{stkl}=D_{st}D_{kl}-D_{kl}D_{st}-(q-1)D_{kt}D_{sl}$, where $s>k$ and $t>l$,
we have two cases to deal with.\par

Case1. If $l\geq j$, then $i>s>k$, $t>l\geq j$, and it follows that
$$\begin{array}{lll}
(f_{ijst},g_{stkl})_{w_1}&=&f_{ijst}D_{kl}-D_{ij}g_{stkl}\\
&=&-qD_{st}D_{ij}D_{kl}+D_{ij}D_{kl}D_{st}+(q-1)D_{ij}D_{kt}D_{sl}\\
&\equiv& -q^2D_{st}D_{kl}D_{ij}+qD_{kl}D_{ij}D_{st}+(q-1)qD_{kt}D_{ij}D_{sl}\\
&\equiv&-q^2D_{kl}D_{st}D_{ij}-q^2(q-1)D_{kt}D_{sl}D_{ij}+q^2(q-1)D_{kt}D_{sl}D_{ij}
+q^2D_{kl}D_{st}D_{ij}\\
&\equiv& 0~\hbox{mod}(S_1,w_1)
\end{array}$$\par

Case2. If $l<j$, then $i>s>k$, $t\ge j>l$, and it follows that
$$\begin{array}{lll}
(f_{ijst},~g_{stkl})_{w_1}&=&f_{ijst}D_{kl}-D_{ij}g_{stkl}\\
&=&-qD_{st}D_{ij}D_{kl}+D_{ij}D_{kl}D_{st}+(q-1)D_{ij}D_{kt}D_{sl}\\
&\equiv&-qD_{st}D_{kl}D_{ij}-q(q-1)D_{st}D_{kj}D_{il}
+D_{kl}D_{ij}D_{st}\\
&{~}&+(q-1)D_{kj}D_{il}D_{st}+(q-1)qD_{kt}D_{ij}D_{sl}\\
&\equiv&-qD_{kl}D_{st}D_{ij}-q(q-1)D_{kt}D_{sl}D_{ij}
-q(q-1)D_{kj}D_{st}D_{il}\\
&{~}&-q(q-1)^2D_{kt}D_{sj}D_{il}+qD_{kl}D_{st}D_{ij}
+q(q-1)D_{kj}D_{st}D_{il}\\
&{~}&+q(q-1)D_{kt}D_{sl}D_{ij}+q(q-1)^2D_{kt}D_{sj}D_{il}\\
&\equiv&0~\hbox{mod}(S_1,w_1)
\end{array}$$\par

$\bullet$ The case $(a\wedge c)$ with $w_1=D_{ij}D_{st}D_{sk}$, where $i>s$, $j\le t,$ $t>k$.\par
Since in this case $w_1=D_{ij}D_{st}D_{sk}=\OV{f}_{ijst}D_{sk}=D_{ij}\OV{h}_{stsk}$, where
$f_{ijst}=D_{ij}D_{st}-qD_{st}D_{ij}$ with $i>s$ and $j\le t$, and $h_{stsk}=D_{st}D_{sk}-D_{sk}D_{st}$ with $t>k$,
we have two cases to deal with.\par

Case 1. If $j\le k$, then $i>s$, $j\le t$, and $j\le k$, thereby
$$\begin{array}{rcl} (f_{ijst},h_{stsk})_{w_1}&=&f_{ijst}D_{sk}-D_{ij}h_{stsk}\\
&=&-qD_{st}D_{ij}D_{sk}+D_{ij}D_{sk}D_{st}\\
&\equiv&-q^2D_{st}D_{sk}D_{ij}+qD_{sj}D_{ij}D_{st}\\
&\equiv&-q^2D_{st}D_{sk}D_{ij}+q^2D_{st}D_{sk}D_{ij}\\
&\equiv&0~\hbox{mod}(S,w_1).\end{array}$$\par

Case 2. If $j>k$, then $i>s$ and $k<j\le t$, thereby
$$\begin{array}{rcl} (f_{ijst}, h_{stsk})_{w_1}&=&f_{ijst}D_{sk}-D_{ij}h_{stsk}\\
&=&-qD_{st}D_{ij}D_{sk}+D_{ij}D_{sk}D_{st}\\
&\equiv&-qD_{st}D_{sk}D_{ij}-q(q-1)D_{st}D_{sj}D_{ik}\\
&{~}&+D_{sk}D_{ij}D_{st}+(q-1)D_{sj}D_{ik}D_{st}\\
&\equiv&-qD_{sk}D_{st}D_{ij}-q(q-1)D_{sj}D_{st}D_{ik}\\
&{~}&+qD_{sk}D_{st}D_{ij}+q(q-1)D_{sj}D_{st}D_{ik}\\
&\equiv&0~\hbox{mod}(S,w_1).\end{array}$$\par

$\bullet$ The case $(b\wedge b)$ with $w=D_{ij}D_{st}D_{kl}$, where $i>s>k,~j>t>l$.\par
Since in this case $w=D_{ij}D_{st}D_{kl}=\OV{g}_{ijst}D_{kl}=D_{ij}\OV{g}_{stkl}$ with
$g_{ijst}=D_{ij}D_{st}-D_{st}D_{ij}-(q-1)D_{sj}D_{it}$ and
$g_{stkl}=D_{st}D_{kl}-D_{kl}D_{st}-(q-1)D_{kt}D_{sl}$, where $i>s>k$ and $j>t>l$,
we have
$$\begin{array}{rcl} (g_{ijst}, g_{stkl})_w&=&g_{ijst}D_{kl}-D_{ij}g_{stkl}\\
&=&-D_{st}D_{ij}D_{kl}-(q-1)D_{sj}D_{it}D_{kl}+D_{ij}D_{kl}D_{st}+(q-1)D_{ij}D_{kt}D_{sl}\\
&\equiv&-D_{st}D_{kl}D_{ij}-(q-1)D_{st}D_{kj}D_{il}-(q-1)D_{sj}D_{kl}D_{it}\\
&{~}&-(q-1)^2D_{sj}D_{kt}D_{il}+D_{kl}D_{ij}D_{st}+(q-1)D_{kj}D_{il}D_{st}\\
&{~}&+(q-1)D_{kt}D_{ij}D_{sl}+(q-1)^2D_{kj}D_{it}D_{sl}\\
&\equiv&-D_{kl}D_{st}D_{ij}-(q-1)D_{kt}D_{sl}D_{ij}-q(q-1)D_{kj}D_{st}D_{il}\\
&{~}&-(q-1)D_{kl}D_{sj}D_{it}-(q-1)^@D_{kj}D_{sl}D_{it}-(q-1)^2D_{kt}D_{sj}D_{il}\\
&{~}&-(q-1)^3D_{kj}D_{st}D_{il}+D_{kl}D_{st}D_{ij}+(q-1)D_{kl}D_{sj}D_{it}\\
&{~}&+q(q-1)D_{kj}D_{st}D_{il}+(q-1)D_{kt}D_{sl}D_{ij}+(q-1)^2D_{kt}D_{sj}D_{il}\\
&{~}&+(q-1)^2D_{kj}D_{sl}D_{it}+(q-1)^3D_{kj}D_{st}D_{il}\\
&\equiv&0~\hbox{mod}(S,w).\end{array}$$\par

$\bullet$ The case $(b\wedge c)$ with $w_1=D_{ij}D_{st}D_{sk}$, where $i>s$, $j>t>k$.\par
Since in this case $w_1=D_{ij}D_{st}D_{sk}=\OV{g}_{ijst}D_{sk}=D_{ij}\OV{h}_{stsk}$, where
$g_{ijst}=D_{ij}D_{st}-D_{st}D_{ij}-(q-1)D_{sj}D_{it}$ with $i>s$ and $j>t$, and
$h_{stsk}=D_{st}D_{sk}-D_{sk}D_{st}$ with $t>k$,
we have
$$\begin{array}{rcl} (g_{ijst}, h_{stsk})_{w_1}&=& g_{ijst}D_{sk}-D_{ij}h_{stsk}\\
&=&-D_{st}D_{ij}D_{sk}-(q-1)D_{sj}D_{it}D_{sk}+D_{ij}D_{sk}D_{st}\\
&\equiv&-D_{st}D_{sk}D_{ij}-(q-1)D_{st}D_{sj}D_{ik}-(q-1)D_{sj}D_{sk}D_{it}\\
&{~}&-(q-1)^2D_{sj}D_{st}D_{ik}+D_{sk}D_{ij}D_{st}+(q-1)D_{sj}D_{ik}D_{st}\\
&\equiv&-D_{sk}D_{st}D_{ij}-(q-1)D_{st}D_{sj}D_{ik}-(q-1)D_{sj}D_{sk}D_{it}\\
&{~}&+(q-1)^2D_{st}D_{sj}D_{ik}+D_{sk}D_{st}D_{ij}+(q-1)D_{sk}D_{sj}D_{it}\\
&{~}&+q(q-1)D_{sj}D_{st}D_{ik}\\
&\equiv&0~\hbox{mod}(S,w_1).\end{array}$$
This finishes the proof of the theorem.\QED\v5

Immediately, Theorem 2.1 gives rise to the following\v5\def\KZ{K\langle D\rangle}

{\bf Corollary 2.2}  The quantized matrix algebra $D_q(n)\cong \KD /J$ has the linear basis, or more precisely, the PBW basis
$$\B =\left\{\left. d^{k_{11}}_{11}d^{k_{12}}_{12}\cdots d^{k_{1n}}_{1n}d^{k_{21}}_{21}\cdots d^{k_{2n}}_{2n}\cdots d^{k_{n1}}_{n1}\cdots d^{k_{nn}}_{nn}~\right |~k_{ij}\in \NZ,(i, j)\in I(n)\right\}.$$ \vskip 6pt

{\bf Proof} With respect to the monomial ordering $\PRC$ on the set $D^*$ of mono words of $K\langle D\rangle$,
we note that
$$\begin{array}{l} D_{11}\PRC D_{12}\PRC\cdots\PRC D_{1n}\PRC
D_{21}\PRC D_{22}\PRC\cdots\PRC D_{2n}\\
\PRC \cdots\PRC
D_{n1}\PRC D_{n2}\PRC\cdots\PRC D_{nn},\end{array}$$
and the Gr\"obner-Shirshov basis $S$ of the ideal $J=\langle S\rangle$ has the set of leading mono words consisting of
$$\begin{array}{l} \OV{f}_{ijst}=D_{ij}D_{st}~\hbox{with}~D_{st}\PRC D_{ij}~\hbox{where} ~s<i,~j\le t,\\
\OV{g}_{ijst}=D_{ij}D_{st}~\hbox{with}~D_{st}\PRC D_{ij}~\hbox{where}~s<i,~t<j,\\
\OV{h}_{ijik}= D_{ij}D_{ik}~\hbox{with}
~D_{ik}\PRC D_{ij},~\hbox{if}~k<j.\end{array}$$
It follows from classical Gr\"obner-Shirshov basis theory that the set of  normal forms of $D^*$ (mod~$S$) is given as follows:
 $$\left\{\left. D^{k_{11}}_{11}D^{k_{12}}_{12}\cdots D^{k_{1n}}_{1n}D^{k_{21}}_{21}\cdots D^{k_{2n}}_{2n}\cdots D^{k_{n1}}_{n1}\cdots D^{k_{nn}}_{nn},~\right |~k_{ij}\in \NZ ,~(i,j)\in I(n)\right\} .$$
Therefore, $D_q(n)$ has the desired PBW basis. \QED\par

\section*{3. $D_q(n)$ is a solvable polynomial algebra}
We start  by recalling from ([9], [12, 15]) the following definitions and notations. Suppose that a finitely generated  $K$-algebra $A=K[a_1,\ldots ,a_n]$ has the PBW $K$-basis $\B =\{ a^{\alpha}=a_{1}^{\alpha_1}\cdots
a_{n}^{\alpha_n}~|~\alpha =(\alpha_1,\ldots ,\alpha_n)\in\NZ^n\}$, and that $\prec$ is a total ordering on $\B$. Then every nonzero element $f\in A$ has a unique expression
$$\begin{array}{rcl} f&=&\lambda_1a^{\alpha (1)}+\lambda_2a^{\alpha (2)}+\cdots +\lambda_ma^{\alpha (m)},\\
&{~}&\hbox{such that}~a^{\alpha (1)}\prec a^{\alpha
(2)}\prec\cdots \prec a^{\alpha (m)},\\
&{~}&\hbox{where}~ \lambda_j\in K^*,~a^{\alpha
(j)}=a_1^{\alpha_{1j}}a_2^{\alpha_{2j}}\cdots a_n^{\alpha_{nj}}\in\B
,~1\le j\le m.
\end{array}$$
Since elements of $\B$ are conventionally called {\it monomials}\index{monomial}, the {\it leading monomial of $f$} is defined as $\LM
(f)=a^{\alpha (m)}$, the {\it leading coefficient of $f$} is defined
as $\LC (f)=\lambda_m$, and the {\it leading term of $f$} is defined
as $\LT (f)=\lambda_ma^{\alpha (m)}$.\v5

{\bf Definition 3.1}  Suppose that the $K$-algebra
$A=K[a_1,\ldots ,a_n]$ has the PBW basis $\B$. If $\prec$ is a
total ordering on $\B$ that satisfies the following three
conditions:{\parindent=1.35truecm\par

\item{(1)} $\prec$ is a well-ordering (i.e., every nonempty subset of $\B$ has a minimal element);\par

\item{(2)} For $a^{\gamma},a^{\alpha},a^{\beta},a^{\eta}\in\B$, if $a^{\gamma}\ne 1$, $a^{\beta}\ne
a^{\gamma}$, and $a^{\gamma}=\LM (a^{\alpha}a^{\beta}a^{\eta})$,
then $a^{\beta}\prec a^{\gamma}$ (thereby $1\prec a^{\gamma}$ for
all $a^{\gamma}\ne 1$);\par

\item{(3)} For $a^{\gamma},a^{\alpha},a^{\beta}, a^{\eta}\in\B$, if
$a^{\alpha}\prec a^{\beta}$, $\LM (a^{\gamma}a^{\alpha}a^{\eta})\ne
0$, and $\LM (a^{\gamma}a^{\beta}a^{\eta})\not\in \{ 0,1\}$, then
$\LM (a^{\gamma}a^{\alpha}a^{\eta})\prec\LM
(a^{\gamma}a^{\beta}a^{\eta})$,\par}{\parindent=0pt
then $\prec$ is called a {\it monomial ordering} on $\B$ (or a
monomial ordering on $A$).} \v5

{\bf Definition 3.2} A finitely generated $K$-algebra $A=K[a_1,\ldots ,a_n]$
is called a {\it solvable polynomial algebra} if $A$ has the PBW $K$-basis $\B =\{
a^{\alpha}=a_1^{\alpha_1}\cdots a_n^{\alpha_n}~|~\alpha
=(\alpha_1,\ldots ,\alpha_n)\in\NZ^n\}$ and a monomial ordering $\prec$ on $\B$, such that
for some $\lambda_{ji}\in K^*$ and  $f_{ji}\in A$,
$$\begin{array}{l} a_ja_i=\lambda_{ji}a_ia_j+f_{ji},~1\le i<j\le n,\\
\LM (f_{ji})\prec a_ia_j~\hbox{whenever}~f_{ji}\ne 0.\end{array}$$\par

Now, we aim to prove the following result.\v5

{\bf Theorem 3.3}  Let $D_q(n)$ be the  quantized matrix algebra over a field $K$, in the sense of [2]. Then $D_q(n)$ is a solvable polynomial algebra in the sense of Definition 3.2. \vskip 6pt

{\bf Proof} Let $I(n)=\{(i,j)~|~i,j=1,2,\cdots,n\}$ with $n\ge 2$. Recall that $D_q(n)$ is the associative $K$-algebra generated by the set of $n^2$ generators
$D=\{ d_{ij}~|~(i,j)\in I(n)\}$  subject to the relations:
$$\begin{array}{ll}
\hbox{F}_1:\quad d_{ij}d_{st}=qd_{st}d_{ij},&\hbox{if}~i>s~\hbox{and}~j\le t,\\
\hbox{F}_2:\quad d_{ij}d_{st}=d_{st}d_{ij}+(q-1)d_{sj}d_{it},&\hbox{if}~i>s~\hbox{and}~j>t,\\
\hbox{F}_3:\quad d_{ij}d_{ik}=d_{ik}d_{ij},&\hbox{for~all}~i,j,k, \end{array}$$
where $i,j,k,s,t=1,2,...,n$ and $q\in K^*$ is the quantum parameter, and
by Corollary 2.2, $D_q(n)$ has the  PBW $K$-basis
$$\B =\left\{\left. d^{k_{11}}_{11}d^{k_{12}}_{12}\cdots d^{k_{1n}}_{1n}d^{k_{21}}_{21}\cdots d^{k_{2n}}_{2n}\cdots d^{k_{n1}}_{n1}\cdots d^{k_{nn}}_{nn}~\right |~k_{ij}\in \NZ,(i, j)\in I(n)\right\}.$$ \par

We now start on constructing a monomial ordering on $\B$ such that all conditions of Definition 3.1 and Definition 3.2 are satisfied. In doing so, we first rewrite $\B$ as
$$\B =\left\{ 1,~d_{i_1j_1}d_{i_2j_2}\cdots d_{i_kj_k}~\left |~
\begin{array}{l} (i_t,j_t)\in I(n),~k\ge 1, \\
(i_1,j_1)\le (i_2,j_2)\le\cdots\le (i_k,j_k)\end{array}\right.\right\} ,$$
where for $(i_{\ell},j_{\ell})$, $(i_t,j_t)\in I(n)$,
$$(i_{\ell},j_{\ell})<(i_t,j_t)\Leftrightarrow\left\{
\begin{array}{l} i_{\ell}<i_t,\\
\hbox{or}~i_{\ell}=i_t~\hbox{and}~j_{\ell}<j_t.\end{array}\right.$$
Then, we define the ordering $\prec$ on the set $D$ of generators: for $d_{kl}$, $d_{ij}\in D$,
$$d_{kl}\prec d_{ij}\Leftrightarrow\left\{\begin{array}{l}
k<i,\\
\hbox{or}~k=i~\hbox{and}~l>j.\end{array}\right.$$
Note that {\it $\prec$ is not the ordering induced by the ordering $<$ on $I(n)$}. Extending the  ordering $\prec$ to $\B$, if we define,
for $u=d_{k_1l_1}d_{k_2l_2}\cdots d_{k_rl_r},$, $v=d_{i_1j_1}d_{i_2j_2}\cdots d_{i_qj_q}\in\B -\{ 1\}$,
$$u\prec v\Leftrightarrow\left\{\begin{array}{l} r<q~\hbox{and}~ d_{k_1l_1}=d_{i_1j_1},d_{k_2l_2}=d_{i_2j_2},\ldots ,d_{k_rl_r}=d_{i_rj_r},\\
\hbox{or there exists an}~m,~1\le m\le r,~\hbox{such that}\\
d_{k_1l_1}=d_{i_1j_1},d_{k_2l_2}=d_{i_2j_2},\ldots ,d_{k_{m-1}l_{m-1}}=d_{i_{m-1}j_{m-1}},~
\hbox{but}~d_{k_ml_m}\prec d_{i_mj_m},\end{array}\right. ,$$
and since conventionally 1 has length 0, we define
$$1\prec u~\hbox{for all}~u=d_{k_1l_1}d_{k_2l_2}\cdots z_{k_rd_r}\in\B -\{ 1\},$$
then it is straightforward to check that $\prec$ is reflexive, antisymmetrical, transitive, and any two elements $u,v\in\B$ are comparable, thereby $\prec$ is a
total ordering on $\B$. Also since $I(n)$ is a finite set, it can be directly verified  that $\prec$ satisfies the descending chain condition on $\B$, namely $\prec$ is a well-ordering on $\B$. \par
It remains to show that $\prec$ satisfies the conditions (2) and (3) of Definition 3.1, and that with respect to $\prec$ on $\B$, the relations $\hbox{F}_1$, $\hbox{F}_2$, and $\hbox{F}_3$  satisfied by generators of $D_q(n)$ have the property required by Definition 3.2. To this end, we first note that Definition 3.2 requires that the product $a_ja_i$ of two generators must be a linear combination of monomials in $\B$, so that $\LM (f_{ji})$ is well defined and $\LM (f_{ji})\prec a_ia_j$ with respect to the ordering $\prec$ defined on $\B$. Bearing in mind this basic requirement, we first observe that in the relations $\hbox{F}_1$, $\hbox{F}_2$, the monomials $d_{st}d_{ij}\in \B$, and in the relation $\hbox{F}_3$, the monomials $d_{ik}d_{ij}\in\B$ for $k<j$. Next, let $d_{ij}, d_{st}, d_{pq}\in D$, and suppose that $d_{st}\prec d_{pq}$. If $(i,j)$ and $(s,t)$ are such that $i>s$ and $j>t$, then the relation $\hbox{F}_2$ gives rise to
$$\begin{array}{rcl} d_{ij}d_{st}&=&d_{st}d_{ij}+(q-1)d_{sj}d_{it}\\
&{~}&\hbox{with}~d_{st}d_{ij}, d_{sj}d_{it}\in\B\\
&{~}&\hbox{and}~ d_{sj}d_{it}\prec d_{st}d_{ij}=\LM (d_{ij}d_{st}).\end{array}$$
On the other hand, if $(i,j)$ and $(p,q)$ are such that $i>p$ and $j>q$, then the relation $\hbox{F}_2$ gives rise to
$$\begin{array}{rcl} d_{ij}d_{pq}&=&d_{pq}d_{ij}+(q-1)d_{pj}d_{iq}\\
&{~}&\hbox{with}~d_{pq}d_{ij}, d_{pj}d_{iq}\in\B\\
&{~}&\hbox{and}~ d_{pj}d_{iq}\prec d_{pq}d_{ij}=\LM (d_{ij}d_{pq}).\end{array}$$
Thus, we have shown that if
$$\begin{array}{l}
(i,j), (s,t)\in I(n)~\hbox{such that}~i>s, j>t,\\
(i,j), (p,q)\in I(n)~\hbox{such that}~i>p, j>q,\end{array}$$
then
$$\begin{array}{l} d_{st}\prec d_{pq}~\hbox{implies}~\LM (d_{ij}d_{st})=
d_{st}d_{ij}\prec d_{pq}d_{ij}=\LM (d_{ij}d_{pq})\\
\hbox{and the generating relations of}~D_q(n)~\hbox{determined by}~\hbox{F}_2\\
\hbox{have the property required by Definition 3.2}.\end{array}\eqno{(1)}$$
Similarly, in the case that
$$\begin{array}{l} (s,t), (i,j)\in I(n)~\hbox{such that}~s>i, t>j,\\
(p,q), (i,j)\in I(n)~\hbox{such that}~ p>i, q>j,\end{array}$$
with the aid of $\hbox{F}_2$ we have
$$\begin{array}{l} d_{st}\prec d_{pq}~\hbox{implies}~\LM (d_{st}d_{ij})=
d_{ij}d_{st}\prec d_{ij}d_{pq}=\LM (d_{pq}d_{ij})\\
\hbox{and the generating relations of}~D_q(n)~\hbox{determined by}~\hbox{F}_2\\
\hbox{have the property required by Definition 3.2}.\end{array}\eqno{(2)}$$
At this stage, bearing in mind the relations $\hbox{F}_1$, $\hbox{F}_2$, $\hbox{F}_3$, and the assertions (1) and (2) derived above, we may conclude that
$$\begin{array}{l} \hbox{for any}~d_{ij}, d_{st}, d_{pq}\in D,~\hbox{if}~d_{st}\prec d_{pq}, ~\hbox{then}\\
\LM (d_{ij}d_{st})\prec \LM (d_{ij}d_{pq}),~\LM (d_{st}d_{ij})\prec\LM (d_{pq}d_{ij}), and\\
\hbox{the generating relations of}~D_q(n)~\hbox{determined by}~\hbox{F}_1, \hbox{F}_2, ~\hbox{and}~\hbox{F}_3\\
\hbox{have the property required by Definition 3.2.}\end{array}\eqno{(3)}$$
Finally, by means of (1), (2), and (3) presented above, it is straightforward to check that the conditions (2) and (3) of Definition 3.1 are satisfied by $\prec$, thereby $\prec$ is a monomial ordering on $\B$, and consequently $D_q(n)$ is a solvable polynomial algebra in the sense of Definition 3.2, as desired.\QED\par

\section*{4. Some applications of the foregoing results}
In this section we specify several applications of the main results of Section 2 and Section 3, so as to demonstrate how  the main results of the foregoing two sections  may bring some perspective of establishing and realizing structural properties of $D_q(n)$ and their modules in a constructive-computational way. All notations used in previous sections are maintained.\v5

We start by recalling three results of [13] in one proposition below, for the reader's convenience.\v5

{\bf Proposition 4.1}  Adopting notations used in [13], let $\KX =K\langle X_1,X_2,\ldots ,X_n\rangle$ be the free $K$-algebra with the set of generators $X=\{ X_1,X_2,\ldots ,X_n\}$, and let $\prec$ be a monomial ordering on $\KX$. Suppose that $\G$ is a Gr\"obner-Shirshov  basis of the ideal $I=\langle \G\rangle$ with respect to $\prec$, such that the set of leading monomials is given by
$$\begin{array}{l} \LM (\G )=\{ X_jX_i~|~1\le i<j\le n\} ,\\
\hbox{or}\\
\LM (\G )=\{ X_iX_j~|~1\le i<j\le n\} .\end{array}$$
Considering the algebra $A=\KX /I$, the following statements hold.\par
(i) [13, P.167, Example 3] The Gelfand-Kirillov dimension GK.dim$A=n$.\par

(ii) [13, P.185, Corollary 7.6] The global homological dimension gl.dim$A=n$, provided $\G$ consists of homogeneous elements with respect to a certain $\NZ$-gradation of $\KX$. (Note that in this case $G^{\NZ}(A)=A$, with the notation used in loc. cit.) \par

(iii) [13, P.201, Corollary 3.2] $A$ is a classical  Koszul algebra,  provided $\G$ consists of quadratic homogeneous elements with respect to the $\NZ$-gradation of $\KX$ such that each $X_i$ is assigned the degree 1, $1\le i\le n$. (Note that in this case $G^{\NZ}(A)=A$, with the notation used in loc. cit.)\par\QED\v5

{\bf Remark} Let $j_1j_2\cdots j_n$ be a permutation of $1, 2, \ldots ,n$. One may notice from the respectively quoted references  in Proposition 4.1 that if, in the case of Proposition 4.1, the monomial ordering $\prec$ employed there is such that
$$\begin{array}{l} X_{j_1}\prec X_{j_2}\prec\cdots \prec X_{j_n},~\hbox{and}\\
\LM (\G )=\{ X_{j_k}X_{j_t}~|~X_{j_t}\prec X_{j_k},~1\le j_k,j_t\le n\} ,\\
\hbox{or}\\
\LM (\G )=\{ X_{j_k}X_{j_t}~|~X_{j_k}\prec X_{j_t},~1\le j_k,j_t\le n\} ,\end{array}$$
then all results still hold true.\v5

Applying Proposition 4.1 and the above remark to $D_q(n)\cong \KZ /J$, we are able to derive  the result below.\v5

{\bf Theorem 4.2} The quantized matrix algebra $D_q(n)$ has the following structural properties.\par
(i) The Hilbert series of $D_q(n)$ is $\frac{1}{(1-t)^{n^2}}$.\par
(ii) The Gelfand-Kirillov dimension GK.dim$D_q(n)=n^2$.\par
(iii) The global homological dimension gl.dim$D_q(n)=n^2$.\par
(iv) $D_q(n)$ is a classical quadratic Koszul algebra.\vskip 6pt

{\bf Proof} Recalling  from Section 2 that with respect to the monomial ordering $\PRC$ on the set $D^*$ of mono words of $\KD$, we have
$$D_{kl}\PRC D_{ij}\Leftrightarrow\left\{
\begin{array}{l} k<i,~l<j,\\
k<i,~j<l,\\
k<i,~j=l,\\
k=i,~l<j.\end{array}\right.\quad (i,j)\in I(n),$$
$$\begin{array}{l} D_{11}\PRC D_{12}\PRC\cdots\PRC D_{1n}\PRC
D_{21}\PRC D_{22}\PRC\cdots\PRC D_{2n}\\
\PRC \cdots\PRC
D_{n1}\PRC D_{n2}\PRC\cdots\PRC D_{nn},\end{array}$$
and thus, the Gr\"obner-Shirshov basis $S$ of the ideal $J=\langle S\rangle$ has the set of leading mono words consisting of
$$\begin{array}{l} \OV{f}_{ijst}=D_{ij}D_{st}~\hbox{with}~D_{st}\PRC D_{ij}~\hbox{where} ~s<i,~j\le t,\\
\OV{g}_{ijst}=D_{ij}D_{st}~\hbox{with}~D_{st}\PRC D_{ij}~\hbox{where}~s<i,~t<j,\\
\OV{h}_{ijik}=D_{ij}D_{ik}~\hbox{with}
~D_{ik}\PRC D_{ij},~\hbox{if}~k<j.\end{array}$$
This means that $D_q(n)$ satisfies the conditions of Proposition 4.1. Therefore, the assertions (i) -- (iv) are established as follows.\par

(i) Since $D_q(n)$ has the PBW $K$-basis as described in Corollary 2.2, it follows that the Hilbert series of $D_q(n)$ is $\frac{1}{(1-t)^{n^2}}$.\par

(ii) This follows from Theorem 2.1, Proposition 4.1(i), and the remark made above.\par

Note that $D_q(n)$ is an $\NZ$-graded algebra defined by a quadratic homogeneous Gr\"obner basis (Theorem 2.1), where each generator $D_{ij}$ is assigned the degree 1, $(i,j)\in I(n)$. The assertions (iii) and (iv) follow from Proposition 4.1(ii) and Proposition 4.1(iii), respectively.\QED\v5

Furthermore, from Theorem 3.3 obtained in the last section we have known that the  quantized matrix algebra $D_q(n)\cong \KZ /J$ is a solvable polynomial algebra in the sense of [9]. Thus, it is well known that every (two-sided, respectively one-sided) ideal of a solvable polynomial algebra $A$ and every submodule of a free (left) $A$-module has a finite Gr\"obner basis with respect to a given monomial ordering, in particular, for one-sided ideals and submodules of free (left) modules there is a noncommutative Buchberger Algorithm which, nowadays, has been successfully implemented in the computer algebra system \textsf{Plural} [11].  At this point, we give several applications of Theorem 3.3 to  $D_q(n)$ and their modules. In what folllows, modules over $D_q(n)$ are meant {\it left $D_q(n)$-modules}.\v5

{\bf Theorem 4.3} Let $D_q(n)$ be the quantized matrix algebra in the sense of [2]. Then the following statements hold.\par

(i) $D_q(n)$ is a Noetherian domain.\par

(ii) Let $L$ be a nonzero left ideal of $D_q(n)$, and $D_q(n)/L$ the left $D_q(n)$-module. Considering Gelfand-Kirillov dimension, we have GK.dim$D_q(n)/L<$ GK.dim$D_q(n)=n^2$, and there is an algorithm for computing GK.dim$D_q(n)$.\par

(iii)  Let $M$ be a finitely generated  $D_q(n)$-module. Then a finite free resolution of $M$ can be algorithmically constructed, and the projective dimension of $M$ can be algorithmically computed. Moreover, every finitely generated projective $D_q(n)$-module $P$ is stably free, thereby the $K_0$-group of $D_q(N)$ is isomorphic to the additive group of integers $\mathbb{Z}$. \par

(iv) Let $M$ be a finitely generated graded $D_q(n)$-module (note that $D_q(n)$ is an $\NZ$-graded algebra in which each generator has degree 1).  Then a minimal homogeneous generating set of $M$ can be algorithmically computed, and a minimal finite graded free resolution of $M$ can be algorithmically constructed.\vskip 6pt

{\bf Proof} (i) Though the property that $D_q(n)$ is a Noetherian domain may be (or may have been) established in some other ways, here we emphasize that this property may also follow immediately from Theorem 3.3. More precisely, the property that $D_q(n)$ has no divisors of zero follows from the fact that $\LM (fg)=\LM(\LM (f)\LM (g))\ne 0$ for all nonzero $f,g\in D_q(n)$, and that the Noetherianess of $D_q(n)$ follows from the fact that every nonzero one-sided ideal of a solvable polynomial algebra has a finite Gr\"obner basis (see [9]). \par

(ii) By Theorem 4.2(ii), GK.dim$D_q(n)=n^2$. Since $D_q(n)$ is a (quadratic) solvable polynomial algebra by Theorem 3.3, it follows from [12, CH.V] that GK.dim$D_q(n)/L<n^2$ (this may also follow from classical Gelfand-Kirillov dimension theory [8], for $D_q(n)$ is now a Noetherian domain), and that there is an algorithm for computing GK.dim$D_q(n)/L$.\par

(iii) This follows from [15, Ch.3].\par

(iv) This follows from [15, Ch.4].\QED\v5

Finally, we show that $D_q(n)$, as a solvable polynomial algebra, also has the elimination property in the sense of [14], and this elimination property can be  realized in a computational way.  To see this clearly, let us first recall the Elimination Lemma given in [14].  Let $A=K[a_1,\ldots ,a_n]$ be a  finitely generated $K$-algebra with
the PBW basis $\B =\{a^{\alpha}=a_1^{\alpha_1}\cdots
a_n^{\alpha_n}~|~\alpha =(\alpha_1,\ldots
,\alpha_n)\in\NZ^n\}$ and, for a subset $U=\{ a_{i_1},...,a_{i_r}\}\subset\{
a_1,...,a_n\}$ with $i_1<i_2<\cdots <i_r$,  let
$$S=\left\{ a_{i_1}^{\alpha_1}\cdots a_{i_r}^{\alpha_r}~\Big |~
(\alpha_1,...,\alpha_r)\in\NZ^r\right\},\quad
V(S)=K\hbox{-span}S.$$

{\bf Lemma 4.4} {[14, Lemma 3.1]} With notation as fixed
above, let $L$ be a nonzero left ideal of $A$ and $A/L$ the left $A$-module defined by $L$. If there is a subset  $U=\{ a_{i_1},\ldots ,a_{i_r}\}
\subset\{a_1,\ldots ,a_n\}$ with $i_1<i_2<\cdots <i_r$, such that $V(S)\cap L=\{ 0\}$, then $$\hbox{GK.dim}(A/L)\ge r.$$  Consequently, if  $A/L$ has finite GK dimension $\hbox{GK.dim}(A/L)=d<n$ ($=$
the number of generators of $A$), then  $$V(S)\cap L\ne
\{ 0\}$$ holds true for every subset $U=\{
a_{i_1},...,a_{i_{d+1}}\}\subset$ $\{ a_1,...,a_n\}$ with
$i_1<i_2<\cdots <i_{d+1}$, in particular, for every $U=\{
a_1,\ldots a_s\}$ with $d+1\le s\le n-1$, we have $V(S)\cap
L\ne \{ 0\}$.\par\QED\v5

For convenience of stating the next theorem, let us write the set of generators of $D_q(n)$ as $D=\{ d_1,d_2,\ldots d_{n^2}\}$, i.e., $D_q(n)=K[d_1,d_2,\ldots ,d_{n^2}]$.\v5

{\bf Theorem 4.5} With notation as fixed above, Let $L$ be a nonzero left ideal of $D_q(n)$. Then the following two statements hold. \par

(i) GK.dim$D_q(n)/L< n^2=$ GK.dim$D_q(n)$. If GK.dim$D_q(n)/L=t$, then
$$V(S)\cap L\ne\{ 0\}$$ holds true for every subset $U=\{
d_{i_1},d_{i_2},...,d_{t+1}\}\subset Z$ with
$i_1<i_2<\cdots <i_{t+1}$, in particular, for every $U=\{
x_1,x_2\ldots x_s\}$ with $t+1\le s\le n^2 -1$, we have $V(S)\cap
L\ne \{ 0\}$.\par

(ii) Without exactly knowing the numerical value GK.dim$D_q(n)/L$, the elimination property for a left ideal  $L=\sum_{i=1}^mD_q(n)\xi_i$ of $D_q(n)$ can be realized in a computational way, as follows:\par

Let $\prec$ be the monomial ordering on the PBW basis $\B$ of $D_q(n)$ as constructed in the proof of Theorem 3.3, and let $V(S)$ be as in (i). Then, employing an elimination ordering  $\lessdot$ with respect to $V(S)$ (which can always  be constructed if the existing monomial ordering on $\B$ is not an elimination ordering, see [15, Proposition 1.6.3]),  a Gr\"obner basis $\G$ of $L$ can be produced by running the noncommutative Buchberger algorithm for solvable polynomial algebras, such that
$$L\cap V(S)\ne \{ 0\} \Leftrightarrow \G\cap V(S)\ne \emptyset .$$\vskip 6pt

{\bf Proof} (i) Since $D_q(n)$ is a solvable polynomial algebra which has the PBW basis $\B$, GK.dim$D_q(n)=n^2$ by Theorem 4.2(ii), and GK.dim$D_q(n)/L<n^2$ by Theorem 4.3(ii), the desired elimination property follows from  Lemma 4.4 mentioned above.\par

(ii) This follows from [15, Corollary 1.6.5].\QED\v5

{\bf Remark} Since $D_q(n)$ is now a solvable polynomial algebra, if  $F=\oplus_{i=1}^sD_q(n)e_i$ is a free (left) $D_q(n)$-module of finite rank, then a similar (even much stronger) result of Theorem 4.6 holds true for any finitely generated submodule $N=\sum_{i=1}^mD_q(n)\xi_i$ of $F$. The reader is referred to [15, Section 2.4] for a detailed argumentation.\v5

\centerline{Reference}{\parindent=.6truecm\par

\item{[1]} L. Bokut et al., {\it Gr\"obner--Shirshov Bases:
Normal Forms, Combinatorial and Decision Problems in Algebra}. World Scientific Publishing,
2020. \url{https://doi.org/10.1142/9287}\par

\item{[2]} R. Dipper and S. Donkin, Quantum $GL_n$. {\it Proc,London Math.Soc.}, 63, 1991,
156--211.

\item{[3]} H. P. Jakobsen and  C. Pagani(2014): Quantized matrix algebras and quantum seeds. {\it Linear and Multilinear Algebra}, DOI: 10.1080/03081087.2014.898297

\item{[4]} H. P. Jakobsen and H. Zhang, The center of the quantized matrix algebra. {\it J. Algebra}, (196)(1997), 458--474.

\item{[5]} H. P. Jakobsen and H. Zhang, The center of the Dipper Donkin quantized matrix algebra. {\it Contributions to Algebra and Geometry}, 2(38), 1997, 411--421.

\item{[6]} H. P. Jakobsen and H. Zhang , A class of quadratic matrix algebras arising from the quantized enveloping algebra $U_q(A_{2n-1})$. {\it J. Math. Phys}., (41)(2000), 2310--2336.

\item{[7]} H. P. Jakobsen, S. J{\o}ndrup, and  A. Jensen, Quadratic algebras of type AIII.III.  In: {\it Tsinghua Science $\&$ Technology}, (3)(1998), 1209--1212 .

\item{[8]} G.R. Krause and T.H. Lenagan, {\it Growth of Algebras and
Gelfand-Kirillov Dimension}. Graduate Studies in Mathematics.
American Mathematical Society, 1991.

\item{[9]} A. Kandri-Rody and V. Weispfenning, Non-commutative
Gr\"obner bases in algebras of solvable type. {\it J. Symbolic
Comput.}, 9(1990), 1--26. Also available as: Technical Report University of Passau,
MIP-8807, March 1988.

\item{[10]} T. Levasseur, Some properties of noncommutative regular graded rings.
{\it Glasgow Math. J}., 34(1992), 277--300.

\item{[11]} V. Levandovskyy and H. Sch\"onemann, Plural: a computer algebra system for noncommutative polynomial algebras. In: {\it Proc. Symbolic and Algebraic Computation}, International Symposium ISSAC 2003, Philadelphia, USA, 176--183, 2003.

\item{[12]} H. Li, {\it Noncommutative Gr\"obner Bases and Filtered-graded Transfer}.
Lecture Notes in Mathematics, Vol. 1795, Springer, 2002.

\item{[13]} H. Li, {\it Gr\"obner Bases in Ring Theory}. World Scientific Publishing Co., 2011. \url{https://doi.org/10.1142/8223}\par

\item{[14]} H. Li, An elimination lemma for algebras with PBW bases. {\it Communications in
Algebra}, 46(8)(2018), 3520--3532.\par

\item{[15]} H. Li, {\it Noncommutative polynomial algebras of solvable type and their modules: Basic constructive-computational theory and methods}. Chapman and Hall/CRC Press, 2021.\par
    
\item{[16]} Li, H., Van Oystaeyen, F. (1996, 2003). {\it Zariskian Filtrations}. K-Monograph in Mathematics, Vol.2. Kluwer Academic Publishers, Berlin Heidelberg: Springer-Verlag.

\end{document}